\documentclass[10pt]{article}
\usepackage{amsfonts,amssymb,amsmath,amsthm,epsfig,euscript,verbatim,graphicx}
\usepackage{textcomp}
\usepackage{enumerate}
\usepackage{epsfig}

\usepackage{tikz}
\usetikzlibrary{patterns}

\usepackage{subcaption}
\usepackage[autostyle]{csquotes}
\usepackage{lmodern}
\usepackage[T1]{fontenc}

\usepackage{palatino}
\usepackage{mathpazo}
\usepackage{inconsolata}
\usepackage{hyperref}

\begingroup
    \makeatletter
    \@for\theoremstyle:=definition,remark,plain\do{%
        \expandafter\g@addto@macro\csname th@\theoremstyle\endcsname{%
            \addtolength\thm@preskip\parskip
            }%
        }
\endgroup

\usepackage[margin=1.1in]{geometry}

\newcommand{\twobyfour}{$2\!\times\!4$}

\makeatletter
\newfont{\footsc}{cmcsc10 at 8truept}
\newfont{\footbf}{cmbx10 at 8truept}
\newfont{\footrm}{cmr10 at 10truept}
\newcommand*\patchAmsMathEnvironmentForLineno[1]{%
  \expandafter\let\csname old#1\expandafter\endcsname\csname #1\endcsname
  \expandafter\let\csname oldend#1\expandafter\endcsname\csname end#1\endcsname
  \renewenvironment{#1}%
     {\linenomath\csname old#1\endcsname}%
     {\csname oldend#1\endcsname\endlinenomath}}%
\newcommand*\patchBothAmsMathEnvironmentsForLineno[1]{%
  \patchAmsMathEnvironmentForLineno{#1}%
  \patchAmsMathEnvironmentForLineno{#1*}}%
\AtBeginDocument{%
\patchBothAmsMathEnvironmentsForLineno{equation}%
\patchBothAmsMathEnvironmentsForLineno{align}%
\patchBothAmsMathEnvironmentsForLineno{flalign}%
\patchBothAmsMathEnvironmentsForLineno{alignat}%
\patchBothAmsMathEnvironmentsForLineno{gather}%
\patchBothAmsMathEnvironmentsForLineno{multline}%
}

\usepackage{mathtools}
\usepackage{parskip}
\usepackage{lineno}
\usepackage[small]{titlesec}

\usepackage[square,comma,numbers,sort&compress]{natbib}

\usepackage{color}
\definecolor{todocolor}{RGB}{205,235,139}
\definecolor{todo-idea}{RGB}{64,150,238}
\definecolor{todo-error}{RGB}{208,31,60}
\definecolor{todo-question}{RGB}{255,255,136}

\usepackage[colorinlistoftodos, color=todocolor, textsize=small]{todonotes}

\newcommand{\Av}{\operatorname{Av}}

\newcommand{\C}{\mathbb{C}}

\newcommand{\CC}{\mathcal{C}}

\renewenvironment{abstract}{
	\begin{list}{}%
	{\setlength{\rightmargin}{1in}%
	\setlength{\leftmargin}{1in}}%
	\item[]\ignorespaces\begin{small}}%
	{\end{small}\unskip\end{list}%
}

\newpagestyle{main}[\small]{
	\headrule
	\sethead[\usepage][][]
	{\sc Completing the Structural Analysis of the 2$\times$4 Permutation Classes}{}{\usepage}
}

\title{\sc Completing the Structural Analysis of the 2$\times$4 Permutation Classes}
\author{
	Samuel Miner\\ \small Mathematics Department\\ \small Pomona College \\ \small Claremont, CA \\ \small \texttt{samuel.miner@gmail.com}
	\and
	Jay Pantone\\ \small Department of Mathematics\\ \small Dartmouth College\\ \small Hanover, NH \\ \small \texttt{jay.pantone@gmail.com}
}

\titleformat{\section}{\large\sc}{\thesection.}{1em}{}
\date{}

\begin{document}
\maketitle

\begin{abstract}
	We study the structure and enumeration of the final two \twobyfour{} permutation classes, completing a research program that has spanned almost two decades. For both classes, careful structural analysis produces a complicated functional equation. One of these equations is solved with the guess-and-check paradigm, while the other is solved with kernel method-like techniques and Gr\"obner basis calculations.
\end{abstract}


\section{Introduction}
\label{section:introduction}

A permutation $\pi = \pi(1)\cdots\pi(n)$ contains a permutation $\sigma = \sigma(1)\cdots\sigma(k)$ if there is a subsequence $\pi(i_1)\cdots\pi(i_k)$ of $\pi$ that is order-isomorphic to $\sigma$---that is, $\pi(i_a) < \pi(i_b)$ if and only if $\sigma(a) < \sigma(b)$. Otherwise, we say that $\pi$ \emph{avoids} $\sigma$. For instance, the permutation $64357218$ contains several occurrences of $231$ but avoids $132$.

Given a set of permutations $B$ such that no permutation in $B$ contains another, define $\Av(B)$ to be the set of all permutations that avoid all of the permutations in $B$. A set of this kind is called a \emph{permutation class}, and the set $B$ is called its \emph{basis}.

The study of permutation classes can be traced back to work of MacMahon~\cite{macmahon:combinatory} in 1915, but began to flourish in the 60s and 70s with the work of Knuth~\cite{knuth:taocp-1}, Pratt~\cite{pratt:computing-permutations}, Tarjan~\cite{tarjan:sorting-networks}, and other on sorting machines. In spite of this long history, much is still unknown about classes avoiding even small sets of permutations. The class $\Av(4231)$ has proved particularly pestiferous; although it has been studied in at least a dozen articles both its generating function and asymptotic behavior remain a mystery. The survey by Vatter~\cite{vatter:perm-survey} provides an easy-to-read comprehensive overview of the field.

A \twobyfour{} permutation class is one whose basis consists of two permutations of length four. There are essentially $56$ different \twobyfour{} classes up to symmetry. These $56$ classes have served as a testing ground for new enumerative techniques for two decades, the first two non-trivial cases appearing in 1998 in works of Atkinson~\cite{atkinson:skew-merged-enum} and B\'ona~\cite{bona:perm-class-smooth}. Of these $56$ classes, there are $38$ different enumerations. In $33$ cases an exact generating function is known, and in $3$ more polynomial-time counting algorithms has been found. In this article, we study the final two \twobyfour{} permutation classes: $\Av(2413,3412)$ and $\Av(3412,4123)$.

%
\section{Preliminaries}
\label{section:preliminaries}

The \emph{diagram} of a permutation is a plot of the points $(i,\pi(i)$ in the Cartesian plane. Often we use partial diagrams, with unshaded boxes that represent where more points may be and shaded boxes that represent areas with no points.

\begin{figure}
\begin{center}
	$\pi\oplus\tau=$
	\begin{tikzpicture}[scale=0.5, baseline=(current bounding box.center)]
		\draw (0,0) rectangle (1,1);
		\draw (1,1) rectangle (2,2);
		\node at (0.5,0.5) {$\pi$};
		\node at (1.5,1.5) {$\tau$};
	\end{tikzpicture}
\quad\quad\quad\quad
	$\pi\ominus\tau=$
	\begin{tikzpicture}[scale=0.5, baseline=(current bounding box.center)]
		\draw (0,1) rectangle (1,2);
		\draw (1,0) rectangle (2,1);
		\node at (0.5,1.5) {$\pi$};
		\node at (1.5,0.5) {$\tau$};
	\end{tikzpicture}
\end{center}
\caption{The sum and skew-sum of permtutations.}
\label{figure:sum-and-skew}
\end{figure}
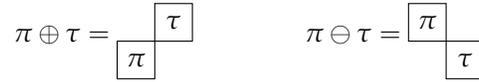

The sum of two permutations $\pi$ and $\tau$ of lengths $k$ and $\ell$ respectively is the permutation $\pi \oplus \ell$ of length $k + \ell$ whose diagram is built by placing the diagram of $\tau$ above and to the right of the diagram of $\pi$, as shown in Figure~\ref{figure:sum-and-skew}. Formally,
\[
    (\pi\oplus\tau)(i)
    =
    \begin{cases}
    \pi(i)&\mbox{for $1\le i\le k$},\\
    \tau(i-k)+k&\mbox{for $k+1\le i\le k+\ell$}.
    \end{cases}
\]
The skew-sum operation $\pi \ominus \tau$, in which the diagram of $\tau$ is placed below and to the right of the diagram of $\pi$ is defined analogously.

A permutation $\pi$ is said to be \emph{sum decomposable} if it can be written as $\pi = \sigma_1 \oplus \sigma_2$ for some permutations $\sigma_1$ and $\sigma_2$ of length at least 1. Otherwise it is \emph{sum indecomposable}. The terms \emph{skew decomposable} and \emph{skew indecomposable} are similarly defined with the skew-sum operation.


\section{$\Av(2413,3412)$}
\label{section:Av(2413,3412)}

We find the generating function for $\CC = \Av(2413, 3412)$ in three steps. First, we describe exactly what permutations in $\CC$ look like, and then we use this description to derive an implicit functional equation for the generating function. Finally, we solve this functional equation with the method of guess-and-check, proving that the generating function for $\CC$ is algebraic of degree $3$ and finding its minimal polynomial.

\subsection{Finding a structural decomposition}
\label{subsection:Av(2413,3412)-structural}

Every permutation $\pi \in \CC$ can be written as the skew sum of skew indecomposable permutations $\pi = \sigma_1 \ominus \cdots \sigma_k$, where $k=1$ when $\pi$ is itself skew indecomposable. At most one of the $\sigma_i$ has length greater than one, or else $\pi$ would contain a $3412$ pattern. Figure~\ref{figure:Av-2413-3412-skew-form} shows the two possible forms of $\pi$: either $\pi$ is strictly decreasing, or $\pi$ contains a single skew indecomposable block preceded and followed by (possibly empty) decreasing permutations. 

Letting $f(z)$ and $f^{\not\ominus}(z)$ denote the generating functions for all permutations and all skew indecomposable permutations in $\CC$, respectively, this proves the identity
\begin{equation}
	\label{equation:f-func-1-no-t}
	f(z) = \frac{1}{1-z} + \frac{f^{\not\ominus}(z)}{(1-z)^2}.
\end{equation}

\begin{figure}
\begin{center}
	\minipage{0.5\textwidth}
	\begin{center}
	\begin{tikzpicture}[scale=1]
			\draw[thick, gray] (0,0) rectangle (2.5,2.5);
			\draw[fill] (0.2, 2.3) circle (2pt);
			\draw[fill] (0.4, 2.1) circle (2pt);
			\draw[fill] (0.6, 1.9) circle (2pt);
			
			\draw[fill] (1.15, 1.35) circle (0.7pt);
			\draw[fill] (1.25, 1.25) circle (0.7pt);
			\draw[fill] (1.35, 1.15) circle (0.7pt);
			
			\draw[fill] (1.9, 0.6) circle (2pt);
			\draw[fill] (2.1, 0.4) circle (2pt);
			\draw[fill] (2.3, 0.2) circle (2pt);
			
			\draw[white,fill] (0.1,4) circle (1pt);
	\end{tikzpicture}
	\qquad
	\begin{tikzpicture}[scale=1]
			\draw[thick, gray] (0,2.5) rectangle (0.9,1.6);
			\draw[fill] (0.15, 2.35) circle (2pt);
			\draw[fill] (0.35, 2.15) circle (0.7pt);
			\draw[fill] (0.45, 2.05) circle (0.7pt);
			\draw[fill] (0.55, 1.95) circle (0.7pt);
			\draw[fill] (0.75, 1.75) circle (2pt);
			
			\begin{scope}[shift={(1.6, -1.6)}]
				\draw[thick, gray] (0,2.5) rectangle (0.9,1.6);
				\draw[fill] (0.15, 2.35) circle (2pt);
				\draw[fill] (0.35, 2.15) circle (0.7pt);
				\draw[fill] (0.45, 2.05) circle (0.7pt);
				\draw[fill] (0.55, 1.95) circle (0.7pt);
				\draw[fill] (0.75, 1.75) circle (2pt);
			\end{scope}
			
			\draw[thick, gray] (0.9,1.6) rectangle (1.6, 0.9) node[midway, black] {$\alpha$};
	\end{tikzpicture}
	\caption{}
	\label{figure:Av-2413-3412-skew-form}
	\end{center}
	\endminipage\hfill
	\minipage{0.5\textwidth}
		\begin{center}
		\begin{tikzpicture}[scale=0.6]
			\draw[thick, lightgray] (3,2) -- (3,8);
			\draw[thick, lightgray] (0,6) -- (8,6);
			\draw[thick, lightgray] (5,2) -- (5,8);
			\fill[pattern=north west lines, pattern color=lightgray] (0,6) rectangle (3,8);
			\fill[pattern=north west lines, pattern color=lightgray] (5,6) rectangle (8,8);
			\fill[pattern=north west lines, pattern color=lightgray] (3,2) rectangle (5,6);
			
			\draw[thick, gray] (0,2) rectangle (8,8);
			
			\draw[fill] (3,8) circle (2pt) node[above] {$\pi(a)$};
			\draw[fill] (1.5,6) circle (2pt) node[above] {$\pi(b)$};
			\draw[fill] (5,7) circle (2pt) node[right] {$\pi(c)$};
			\draw[thick] (5.2, 5.8) -- (7.8, 2.2);
			
			\node at (4, 7) {$R_1$};
			\node at (1.5, 4) {$R_2$};
			\node[above right] at (6.5, 4) {$R_3$};
			
		\end{tikzpicture}
		\caption{}
		\label{figure:Av-2413-3412-skew-1}
	\end{center}
	\endminipage
\end{center}
\end{figure}

Now we begin the much more laborious task of describing the skew indecomposable permutations in $\CC$. Suppose $\pi = \pi(1)\cdots\pi(n)$ is such a permutation. Let $\pi(a)$ be the maximum entry in the permutation, noting that $a \neq 1$, or else $\pi$ would begin with its largest element and hence be skew decomposable. Therefore we can define $\pi(b)$ to be the entry of largest value to the left of $\pi(a)$. Similarly, let $\pi(c)$ be the rightmost element that is both to the right of $\pi(a)$ and above $\pi(b)$ (if such an element exists).

Figure~\ref{figure:Av-2413-3412-skew-1} shows what the plot of such a permutation looks like. The shaded regions indicate areas where there cannot be entries of $\pi$. The area that lies to the left of $\pi(a)$ and above $\pi(b)$ is shaded by our choice of $\pi(b)$. The area to the right of $\pi(c)$ and above $\pi(b)$ is shaded by our choice of $\pi(c)$. Lastly, the area that lies below $\pi(b)$ and horizontally between $\pi(a)$ and $\pi(c)$ is shaded because any entry here would play the role of a $1$ in a $2413$ pattern (together with $\pi(b)$, $\pi(a)$, and $\pi(c)$).

The avoidance of the patterns $2413$ and $3412$ allows us to draw conclusions about the subpermutations that can lie in each of the three regions $R_1$, $R_2$, and $R_3$ shown in Figure~\ref{figure:Av-2413-3412-skew-1}. The subpermutation in $R_1$, together with the entry $\pi(c)$ can be any permutation in $\CC$. The same is true for the subpermutation in $R_2$ together with the entry $\pi(b)$. On the other hand, the diagonal line in the region $R_3$ indicates that the subpermutation consisting of the entries in this cell must be decreasing---otherwise a $3412$ pattern is created.

There are two cases to consider: either $R_3$ is empty or nonempty. In the case where $R_3$ is empty, we have completely described the structure of $\pi$. More precisely, in this case $\pi = \sigma \oplus (1 \ominus \tau)$, for some permutations $\sigma,\tau \in \CC$ such that $\sigma$ is nonempty. Having dispensed of this first case, now assume that $R_3$ is nonempty, and let $\pi(d)$ be the topmost entry in $R_3$. Figure~\ref{figure:Av-2413-3412-skew-2} shows this configuration.

\begin{figure}
	\begin{center}
	\minipage{0.5\textwidth}
		\begin{center}
		\begin{tikzpicture}[scale=0.8]
			\draw[thick, lightgray] (3,2) -- (3,8);
			\draw[thick, lightgray] (0,6) -- (8,6);
			\draw[thick, lightgray] (5,2) -- (5,8);
			\draw[thick, lightgray] (0,5) -- (8,5);
			\fill[pattern=north west lines, pattern color=lightgray] (0,6) rectangle (3,8);
			\fill[pattern=north west lines, pattern color=lightgray] (5,6) rectangle (8,8);
			\fill[pattern=north west lines, pattern color=lightgray] (3,2) rectangle (5,6);
			\fill[pattern=north west lines, pattern color=lightgray] (5,5) rectangle (8,6);
			
			\draw[thick, gray] (0,2) rectangle (8,8);
			
			\draw[fill] (3,8) circle (2pt) node[above] {$\pi(a)$};
			\draw[fill] (1.5,6) circle (2pt) node[above] {$\pi(b)$};
			\draw[fill] (5,7) circle (2pt) node[right] {$\pi(c)$};
			\draw[fill] (5.4, 5) circle (2pt) node[right=4pt, above] {$\pi(d)$};	
			\draw[thick] (5.4, 5) -- (7.8, 2.2);
			
			\node at (4, 7) {$R_1$};
			\node at (1.5, 5.5) {$R_2$};
			\node[above right] at (6.5, 4) {$R_3$};
			\node at (1.5, 3.5) {$R_4$};
		
		\end{tikzpicture}
		\caption{}
		\label{figure:Av-2413-3412-skew-2}
		\end{center}
	\endminipage\hfill
	\minipage{0.5\textwidth}
	\begin{center}
		\begin{tikzpicture}[scale=0.8]
			\draw[thick, lightgray] (3,2) -- (3,8);
			\draw[thick, lightgray] (0,7) -- (8,7);
			\draw[thick, lightgray] (5,2) -- (5,8);
			\draw[thick, lightgray] (0,5) -- (8,5);
			\draw[thick, lightgray] (2,2) -- (2, 5);
			\draw[thick, lightgray, dashed] (2,5.1) -- (2, 6.5);
			\draw[thick, lightgray, dashed] (0.4,2)--(0.4,6.75);
			\draw[thick, lightgray, dashed] (0.8,2)--(0.8,6.25);
			\draw[thick, lightgray, dashed] (1.2,2)--(1.2,6);
			\draw[thick, lightgray, dashed] (1.6,2)--(1.6,5.5);
			\draw[thick, lightgray, dashed] (2,2)--(2,6.5);
			\draw[thick, lightgray, dashed] (0, {5-2*2.8/2.4})--(7.4,{5-2*2.8/2.4});

			\fill[pattern=north west lines, pattern color=lightgray] (0,7) rectangle (3,8);
			\fill[pattern=north west lines, pattern color=lightgray] (5,7) rectangle (8,8);
			\fill[pattern=north west lines, pattern color=lightgray] (3,2) rectangle (5,7);
			\fill[pattern=north west lines, pattern color=lightgray] (5,5) rectangle (8,7);
			\fill[pattern=north west lines, pattern color=lightgray] (2,2) rectangle (3,5);
			
			\draw[thick, gray] (0,2) rectangle (8,8);
			
			\draw[fill] (3,8) circle (2pt) node[above] {$\pi(a)$};
			\draw[fill] (2.5,7) circle (2pt) node[above] {$\pi(b)$};
			\draw[fill] (5,7.5) circle (2pt) node[right] {$\pi(c)$};
			\draw[fill] (5.4, 5) circle (2pt) node[right=4pt, above] {$\pi(d)$};	
			\draw[fill] (7.4, {5-2*2.8/2.4}) circle (2pt) node[above left=0pt and 10pt] {$\pi(e)$};	
			
			\draw[thick] (5.4, 5) -- (7.8, 2.2);
			
			\draw[fill] (0.4, 6.75) circle (2pt);
			\draw[fill] (0.8, 6.25) circle (2pt);
			\draw[fill] (1.2, 6) circle (2pt);
			\draw[fill] (1.6, 5.5) circle (2pt);
			\draw[fill] (2, 6.5) circle (2pt);
			
			\node at (4, 7.5) {$R_1$};
			\node[right=3pt] at (2, 6) {$R_2$};
			\node[above right] at (6.8, 3.8) {$R_3$};
			\node at (0.2, 4.8) {\footnotesize $R_4$};
			
		\end{tikzpicture}
		\caption{}
		\label{figure:Av-2413-3412-skew-3}
	\end{center}
	\endminipage	
\end{center}
\end{figure}

While the conditions on $R_1$, $R_2$, and $R_3$ remain the same, the set of subpermutations that can be drawn on $R_4$ is highly dependent on the subpermutations that appear in $R_2$ and $R_3$. The subpermutation in $R_2$ does not have to be decreasing, but will have some initial run of decreasing entries of length at least one before the leftmost ascent.  As demonstrated in Figure~\ref{figure:Av-2413-3412-skew-3}, there can be no entries in $R_4$ lying to the right of the leftmost ascent in $R_2$; otherwise a $3412$ pattern is created.

Our next goal is to specify what permutations can lie in the region $R_4$ vertically between $\pi(d)$, and the second-highest entry in the region $R_3$ which we've labeled $\pi(e)$ in Figure~\ref{figure:Av-2413-3412-skew-3}. (Alternatively, if there is no second entry in $R_3$, then we are specifying the permutations that lie vertically below $\pi(d)$.)

The example permutation shown in Figure~\ref{figure:Av-2413-3412-skew-3} has an initial decreasing run of length four in the region $R_2$, inducing five vertical strips in the region $R_4$. First note that if there is a decreasing pair of entries in the region $R_4$ that crosses over the dividing line between two vertical strips, then a $2413$ pattern is created---the decreasing pair constitutes the $2$ and the $1$, the entry in $R_2$ whose vertical line they cross plays the role of $4$, and $\pi(d)$ plays the role of $3$. Hence, all the entries in each vertical strip of $R_4$ must lie both above and to the right of all entries in the strips to its left.

Figure~\ref{figure:Av-2413-3412-skew-4} shows where entries can occur in the region $R_4$ for the sample permutation from Figure~\ref{figure:Av-2413-3412-skew-3}. Each of these five boxes can contain any permutation in the class $\CC$ (including the empty permutation). It's important to note that there are five boxes precisely because the subpermutation of all entries above $\pi(d)$ has an initial decreasing run of length four.

\begin{figure}
\begin{center}
	\begin{tikzpicture}[scale=0.8]
		\draw[thick, lightgray] (3,2) -- (3,8);
		\draw[thick, lightgray] (0,7) -- (8,7);
		\draw[thick, lightgray] (5,2) -- (5,8);
		\draw[thick, lightgray] (0,5) -- (8,5);
		\draw[thick, lightgray] (2,2) -- (2, 5);
		\draw[thick, lightgray, dashed] (2,5.1) -- (2, 6.5);
		\draw[thick, lightgray, dashed] (0.4,2)--(0.4,6.75);
		\draw[thick, lightgray, dashed] (0.8,2)--(0.8,6.25);
		\draw[thick, lightgray, dashed] (1.2,2)--(1.2,6);
		\draw[thick, lightgray, dashed] (1.6,2)--(1.6,5.5);
		\draw[thick, lightgray, dashed] (2,2)--(2,6.5);
		\draw[thick, lightgray, dashed] (0, {5-2*2.8/2.4})--(7.4,{5-2*2.8/2.4});

		\fill[pattern=north west lines, pattern color=lightgray] (0,7) rectangle (3,8);
		\fill[pattern=north west lines, pattern color=lightgray] (5,7) rectangle (8,8);
		\fill[pattern=north west lines, pattern color=lightgray] (3,2) rectangle (5,7);
		\fill[pattern=north west lines, pattern color=lightgray] (5,5) rectangle (8,7);
		\fill[pattern=north west lines, pattern color=lightgray] (2,2) rectangle (3,5);
		
		\fill[pattern=north west lines, pattern color=lightgray] (0,{5-2*2.8/2.4}) rectangle (2,5); 
		
		\draw[thick, gray, fill=white] (0,{5-2*2.8/2.4}) rectangle (0.4,{5-2*2.8/2.4+0.4});
		\draw[thick, gray, fill=white] (0.4,{5-2*2.8/2.4+0.4}) rectangle (0.8,{5-2*2.8/2.4+0.8});
		\draw[thick, gray, fill=white] (0.8,{5-2*2.8/2.4+0.8}) rectangle (1.2,{5-2*2.8/2.4+1.2});
		\draw[thick, gray, fill=white] (1.2,{5-2*2.8/2.4+1.2}) rectangle (1.6,{5-2*2.8/2.4+1.6});
		\draw[thick, gray, fill=white] (1.6,{5-2*2.8/2.4+1.6}) rectangle (2,{5-2*2.8/2.4+2});
		
		\draw[thick, gray] (0,2) rectangle (8,8);
		
		\draw[fill] (3,8) circle (2pt) node[above] {$\pi(a)$};
		\draw[fill] (2.5,7) circle (2pt) node[above] {$\pi(b)$};
		\draw[fill] (5,7.5) circle (2pt) node[right] {$\pi(c)$};
		\draw[fill] (5.4, 5) circle (2pt) node[right=4pt, above] {$\pi(d)$};	
		\draw[fill] (7.4, {5-2*2.8/2.4}) circle (2pt) node[above left=0pt and 10pt] {$\pi(e)$};	
		
		\draw[thick] (5.4, 5) -- (7.8, 2.2);
		
		\draw[fill] (0.4, 6.75) circle (2pt);
		\draw[fill] (0.8, 6.25) circle (2pt);
		\draw[fill] (1.2, 6) circle (2pt);
		\draw[fill] (1.6, 5.5) circle (2pt);
		\draw[fill] (2, 6.5) circle (2pt);
		
		\node at (4, 7.5) {$R_1$};
		\node[right=3pt] at (2, 6) {$R_2$};
		\node[above right] at (6.8, 3.8) {$R_3$};
		\node at (0.2, 4.8) {\footnotesize $R_4$};
		
	\end{tikzpicture}
	\caption{}
	\label{figure:Av-2413-3412-skew-4}
\end{center}
\end{figure}

Having categorized the subpermutations that may lie in the region $R_4$ vertically between $\pi(d)$ and $\pi(e)$ (if it exists), we now move on to consider the entries in the region $R_4$ that lie between $\pi(e)$ (if it exists) and the next lower entry in $R_3$. These entries obey the same constraints as the other entries in $R_4$, replacing $\pi(d)$ by $\pi(e)$: they cannot lie to the right of the leftmost ascent of the subpermutation composed of all points lying above $\pi(e)$. Note that the initial decreasing run of this subpermutation may contain some points from the initial decreasing run of the subpermutation lying above $\pi(d)$ together with some points from the leftmost nonempty box in $R_4$ lying vertically between $\pi(d)$ and $\pi(e)$.  

This same structure of entries in $R_4$ is repeated between every pair of entries in $R_3$, and once more below the bottommost entry of $R_3$. We use the term \emph{layer} to mean a nonempty group of entries in $R_4$ lying between a pair of $R_3$ entries together with the largest nonempty sequence of consecutive entries in $R_3$ directly above these entries. Figure~\ref{figure:Av-2413-3412-skew-5} shows an example with three layers. 

\begin{figure}
\begin{center}
	\begin{tikzpicture}[scale=1]
		\draw[thick, lightgray] (3,0) -- (3,8);
		\draw[thick, lightgray] (-2,7) -- (8,7);
		\draw[thick, lightgray] (5,0) -- (5,8);
		\draw[thick, lightgray] (-2,5) -- (8,5);
		\draw[thick, lightgray] (2,0) -- (2, 5);
		\draw[thick, lightgray, dashed] (2,5.1) -- (2, 6.5);
		\draw[thick, lightgray, dashed] (0.4-2,0)--(0.4-2,6.75);
		\draw[thick, lightgray, dashed] (0.8-2,0)--(0.8-2,6.25);
		\draw[thick, lightgray, dashed] (1.2,0)--(1.2,6);
		\draw[thick, lightgray, dashed] (1.6,0)--(1.6,5.5);
		\draw[thick, lightgray, dashed] (2,0)--(2,6.5);
		\draw[thick, lightgray, dashed] (-2, {5-1.1*4.8/2.4})--(6.57,{5-1.1*4.8/2.4});
		\draw[thick, lightgray, dashed] (-2, {5-2*4.8/2.4})--(7.4,{5-2*4.8/2.4});
		\draw[thick, lightgray, dashed] (-2, {5-(1.3)*4.8/2.4}) -- (6.7, {5-(1.3)*4.8/2.4});

		\fill[pattern=north west lines, pattern color=lightgray] (-2,7) rectangle (3,8);
		\fill[pattern=north west lines, pattern color=lightgray] (5,7) rectangle (8,8);
		\fill[pattern=north west lines, pattern color=lightgray] (3,0) rectangle (5,7);
		\fill[pattern=north west lines, pattern color=lightgray] (5,5) rectangle (8,7);
		\fill[pattern=north west lines, pattern color=lightgray] (2,0) rectangle (3,5);
		\fill[pattern=north west lines, pattern color=lightgray] (-2,{5-2*2.8/2.4}) rectangle (2,5); 
		\fill[pattern=north west lines, pattern color=lightgray] (-1.2,0) rectangle (3,{5-2*4.8/2.4});
		\fill[pattern=north west lines, pattern color=lightgray] (-2,0) rectangle (-1.2,{5-(2)*4.8/2.4});
		
		\draw[thick, gray, pattern=north west lines, pattern color=lightgray] (0-2,{5-1.1*4.8/2.4}) rectangle (0.4-2,{5-1.1*4.8/2.4+0.33});
		\draw[thick, gray, pattern=north west lines, pattern color=lightgray] (0.4-2,{5-1.1*4.8/2.4+0.33}) rectangle (0.8-2,{5-2*2.8/2.4+0.8});
		\draw[thick, gray, fill=white] (0.8-2,{5-2*2.8/2.4+0.8}) rectangle (1.2,{5-2*2.8/2.4+1.2});
		\draw[thick, lightgray, dashed] (0.92-1.4, 0) -- (0.92-1.4, {5-2*2.8/2.4+1.08});
		\draw[thick, lightgray, dashed] (1.08-0.6, 0) -- (1.08-0.6, {5-2*2.8/2.4+0.92});
		\draw[thick, gray] (0.8-2,{5-2*2.8/2.4+0.8}) rectangle (1.2,{5-2*2.8/2.4+1.2});
		\draw[thick, gray, fill=white] (1.2,{5-2*2.8/2.4+1.2}) rectangle (1.6,{5-2*2.8/2.4+1.6});
		\draw[thick, gray, fill=white] (1.6,{5-2*2.8/2.4+1.6}) rectangle (2,{5-2*2.8/2.4+2});
		
		\draw[thick, gray, fill=white] (-2,{5-2.3*4.8/2.4}) rectangle (-1.6,{5-2.3*4.8/2.4+0.2});
		\draw[thick, gray, fill=white] (-1.6,{5-2.3*4.8/2.4+0.2}) rectangle (-1.4,{5-2.3*4.8/2.4+0.4});
		\draw[thick, gray, fill=white] (-1.4,{5-2.3*4.8/2.4+0.4}) rectangle (-1.2,{5-2.3*4.8/2.4+0.6});
		
		\fill[pattern=north west lines, pattern color=lightgray] (-2,{5-2*4.8/2.4}) rectangle (3,{5-2*4.8/2.4+2});
		
		\draw[thick, gray, pattern=north west lines, pattern color=lightgray] (-2,{5-2*4.8/2.4}) rectangle (-1.6,{5-2*4.8/2.4+0.3});
		\draw[thick, gray, fill=white] (-1.6,{5-2*4.8/2.4+0.3}) rectangle (-1.2,{5-2*4.8/2.4+0.6});
		\draw[thick, lightgray, dashed] (-1.4, 0) -- (-1.4, {5-2*4.8/2.4+0.45});
		\draw[thick, gray] (-1.6,{5-2*4.8/2.4+0.3}) rectangle (-1.2,{5-2*4.8/2.4+0.6});
		\draw[thick, gray, fill=white] (-1.2,{5-2*4.8/2.4+0.6}) rectangle (0.92-1.4,{5-2*4.8/2.4+0.9});
		\draw[thick, gray, fill=white] (0.92-1.4,{5-2*4.8/2.4+0.9}) rectangle (1.08-0.6,{5-2*4.8/2.4+1.2});
		\draw[thick, gray, pattern=north west lines, pattern color=lightgray] (1.08-0.6,{5-2*4.8/2.4+1.2}) rectangle (1.2,{5-1.3*4.8/2.4});

		\draw[fill] (0.92-1.4, {5-2*2.8/2.4+1.08}) circle (1.5pt);
		\draw[fill] (1.08-0.6, {5-2*2.8/2.4+0.92}) circle (1.5pt);
		\draw[fill] (1.4, {5-2*2.8/2.4+1.4}) circle (1.5pt);
		\draw[fill] (1.72, {5-2*2.8/2.4+1.72}) circle (1.5pt);
		\draw[fill] (1.88, {5-2*2.8/2.4+1.88}) circle (1.5pt);
		\draw[fill] (-1.4, {5-2*4.8/2.4+0.45}) circle (1.5pt);
		\draw[fill] (-0.84, {5-2*4.8/2.4+0.75}) circle (1.5pt);
		\draw[fill] (0, {5-2*4.8/2.4+1.05}) circle (1.5pt);

		\draw[thick, gray] (-2,0) rectangle (8,8);
		
		\draw[fill] (3,8) circle (1.5pt) node[above] {$\pi(a)$};
		\draw[fill] (2.5,7) circle (1.5pt) node[above] {$\pi(b)$};
		\draw[fill] (5,7.5) circle (1.5pt) node[right] {$\pi(c)$};
		\draw[fill] (5.4, 5) circle (1.5pt) node[above] {$\pi(d)$};	
		\draw[fill] (6.5, {5-1.1*4.8/2.4}) circle (1.5pt) node[above right] {$\pi(e)$};
		\draw[fill] (6.6, {5-1.2*4.8/2.4}) circle (1.5pt);
		\draw[fill] (6.7, {5-1.3*4.8/2.4}) circle (1.5pt);
		\draw[fill] (7.4, {5-2*4.8/2.4}) circle (1.5pt);

		\draw[thick, gray] (8.2, 5.1) -- (8.2, 3.1);
		\draw[thick, gray] (8.05, 5.1) -- (8.35, 5.1);
		\draw[thick, gray] (8.05, 3.1) -- (8.35, 3.1);
		\node[rotate=270, gray] at (8.4, 4.1) {\texttt{layer one}};
		
		\draw[thick, gray] (8.2, 3.0) -- (8.2, 1.2);
		\draw[thick, gray] (8.05, 3.0) -- (8.35, 3.0);
		\draw[thick, gray] (8.05, 1.2) -- (8.35, 1.2);
		\node[rotate=270, gray] at (8.4, 2.1) {\texttt{layer two}};
		
		\draw[thick, gray] (8.2, 1.1) -- (8.2, 0.1);
		\draw[thick, gray] (8.05, 1.1) -- (8.35, 1.1);
		\draw[thick, gray] (8.05, 0.1) -- (8.35, 0.1);
		\node[rotate=270, gray] at (8.65, 0.625) {\texttt{layer}};
		\node[rotate=270, gray] at (8.4, 0.625) {\texttt{three}};

		\draw[fill] (-1.8, {5-2.3*4.8/2.4+0.1}) circle (1.5pt);
		\draw[fill] (-1.5, {5-2.3*4.8/2.4+0.3}) circle (1.5pt);
		\draw[fill] (-1.3, {5-2.3*4.8/2.4+0.5}) circle (1.5pt);

		\draw[thick, gray] (-1.6,{5-2.3*4.8/2.4+0.2}) rectangle (-1.4,{5-2.3*4.8/2.4+0.4});
		\draw[thick, gray] (-1.4,{5-2.3*4.8/2.4+0.4}) rectangle (-1.2,{5-2.3*4.8/2.4+0.6});

		\draw[thick] (5.4, 5) -- (7.8, 0.2);
		
		\draw[fill] (0.4-2, 6.75) circle (1.5pt);
		\draw[fill] (0.8-2, 6.25) circle (1.5pt);
		\draw[fill] (1.2, 6) circle (1.5pt);
		\draw[fill] (1.6, 5.5) circle (1.5pt);
		\draw[fill] (2, 6.5) circle (1.5pt);
		
		\node at (4, 7.5) {$R_1$};
		\node at (2.5, 6) {$R_2$};
		\node[above right] at (6.8, 3.8) {$R_3$};
		\node at (0.2-2, 4.8) {$R_4$};
		
	\end{tikzpicture}
	\caption{}
	\label{figure:Av-2413-3412-skew-5}
\end{center}
\end{figure}

\subsection{Constructing a functional equation}

The beginning of Section~\ref{subsection:Av(2413,3412)-structural} makes it clear that if the enumeration for the skew indecomposable permutations can be found, then the enumeration for the skew decomposable permutations (and therefore that of the whole class) will quickly follow. We will construct a functional equation for the skew indecomposable permutations by observing that every such permutation either:
\begin{enumerate}[(i)]
	\item has no entries in region $R_3$ of Figure~\ref{figure:Av-2413-3412-skew-1} (i.e., has no layers), or
	\item has at least one entry in region $R_3$ of Figure~\ref{figure:Av-2413-3412-skew-1} (i.e., has at least one layer), and thus can be built by adding a layer to a skew indecomposable permutation with one fewer layer.
\end{enumerate}

In order to add a layer to a skew indecomposable permutation, we must know the length of the initial decreasing run. To this end, consider the bivariate generating functions $f(z,t)$ and $f^{\not\ominus}(z,t)$ that count permutations in $\CC$ and skew indecomposable permutations in $\CC$, respectively, where $z$ tracks the length of the permutation and $t$ tracks the length of the initial decreasing run. Recall that every permutation has an initial decreasing run of length at least one, and the number of vertical strips (in the sense of Figure~\ref{figure:Av-2413-3412-skew-4}) is equal to one more than the length of the initial decreasing run.

The decomposition into cases (i) and (ii) above translates to the equation
\[
	f^{\not\ominus}(z,t) = g_{(i)}(z,t) + g_{(ii)}(z,t),
\]
where $g_{(i)}(z,t)$ and $g_{(ii)}(z,t)$ count the permutations in cases (i) and (ii) respectively. 

It was previously stated that the permutations in case (i) can be written in the form $\pi = \sigma \oplus (1 \ominus \tau)$ where $\sigma,\tau \in \CC$ and $\sigma$ is nonempty. Note that the initial decreasing run of $\pi$ is composed exactly of the initial decreasing run of $\sigma$, and that the initial decreasing run of $\tau$ is irrelevant. Therefore,
\[
	g_{(i)}(z,t) = \underbrace{(f(z,t)-1)}_{\sigma}\cdot \underbrace{\left(z \cdot f(z,1)\right)}_{1 \ominus \tau} = zf(z,1)(f(z,t)-1).
\]

Of course, $g_{(ii)}(z,t)$ is more complicated. Let us first consider a particular example. Suppose that $\pi$ is a skew indecomposable permutation in $\CC$ with an initial decreasing run with three entries. To add a layer to $\pi$ we first add one or more decreasing entries in the region $R_3$, and then we fill in the four vertical strips created by the three initial decreasing entries. 

Of the four subpermutations (which we'll now call the four \emph{components}) inserted into each of the four vertical strips, at least one must be nonempty. The new initial decreasing run of the permutation resulting from adding a layer to $\pi$ may consist of some entries from the initial decreasing run of $\pi$ as well as some of the new layer. Figure~\ref{figure:Av-2413-3412-skew-new-layer} shows, for this example, the four possible cases, depending on which component is the leftmost nonempty component. The presence of a white box with a point indicates that the component is nonempty.

\begin{figure}
\begin{center}
	\begin{tikzpicture}[scale=0.9]
		
		\fill[thick, gray, pattern=north west lines, pattern color=lightgray] (0.7, 0.8) rectangle (3.4, 3.6);
		
		\draw[thick, gray, dashed] (1.2, 5) -- (1.2, 0.8);
		\draw[thick, gray, dashed] (1.7, 4.2) -- (1.7, 0.8);
		\draw[thick, gray, dashed] (2.2, 3.8) -- (2.2, 0.8);
		\draw[thick, gray, dashed] (2.7, 4.6) -- (2.7, 0.8);

		\draw[fill] (1.2,5) circle (2pt);
		\draw[fill] (1.7,4.2) circle (2pt);
		\draw[fill] (2.2,3.8) circle (2pt);
		\draw[fill] (2.7,4.6) circle (2pt);
		\draw[fill] (3.6,3.4) circle (2pt);
		\draw[fill] (4.0,3) circle (2pt);
		\draw[thick, gray] (0.7, 3.6) rectangle (3.4, 5.2);
		\draw[thick, gray] (3.4,3.6) rectangle (4.2, 2.8);

		\draw[thick, gray, fill=white] (0.7, 0.8) rectangle (1.2, 1.3);
		\draw[thick, gray, fill=white] (1.2, 1.3) rectangle (1.7, 1.8);
		\draw[thick, gray, fill=white] (1.7, 1.8) rectangle (2.2, 2.3);
		\draw[thick, gray, fill=white] (2.2, 2.3) rectangle (2.7, 2.8);
		
		\draw[fill] (0.95, 1.05) circle (2pt);
		\draw[fill] (1.45, 1.55) circle (2pt);
		\draw[fill] (1.95, 2.05) circle (2pt);
		\draw[fill] (2.45, 2.55) circle (2pt);

	\end{tikzpicture}
	\quad
	\begin{tikzpicture}[scale=0.9]
		
		\fill[thick, gray, pattern=north west lines, pattern color=lightgray] (0.7, 0.8) rectangle (3.4, 3.6);
		
		\draw[thick, gray, dashed] (1.2, 5) -- (1.2, 0.8);
		\draw[thick, gray, dashed] (1.7, 4.2) -- (1.7, 0.8);
		\draw[thick, gray, dashed] (2.2, 3.8) -- (2.2, 0.8);
		\draw[thick, gray, dashed] (2.7, 4.6) -- (2.7, 0.8);

		\draw[fill] (1.2,5) circle (2pt);
		\draw[fill] (1.7,4.2) circle (2pt);
		\draw[fill] (2.2,3.8) circle (2pt);
		\draw[fill] (2.7,4.6) circle (2pt);
		\draw[fill] (3.6,3.4) circle (2pt);
		\draw[fill] (4.0,3) circle (2pt);
		\draw[thick, gray] (0.7, 3.6) rectangle (3.4, 5.2);
		\draw[thick, gray] (3.4,3.6) rectangle (4.2, 2.8);
		
		\draw[thick, gray, fill=white] (1.2, 1.3) rectangle (1.7, 1.8);
		\draw[thick, gray, fill=white] (1.7, 1.8) rectangle (2.2, 2.3);
		\draw[thick, gray, fill=white] (2.2, 2.3) rectangle (2.7, 2.8);
		
		\draw[fill] (1.45, 1.55) circle (2pt);
		\draw[fill] (1.95, 2.05) circle (2pt);
		\draw[fill] (2.45, 2.55) circle (2pt);

	\end{tikzpicture}
	\quad
	\begin{tikzpicture}[scale=0.9]
		
		\fill[thick, gray, pattern=north west lines, pattern color=lightgray] (0.7, 0.8) rectangle (3.4, 3.6);
		
		\draw[thick, gray, dashed] (1.2, 5) -- (1.2, 0.8);
		\draw[thick, gray, dashed] (1.7, 4.2) -- (1.7, 0.8);
		\draw[thick, gray, dashed] (2.2, 3.8) -- (2.2, 0.8);
		\draw[thick, gray, dashed] (2.7, 4.6) -- (2.7, 0.8);
				
		\draw[fill] (1.2,5) circle (2pt);
		\draw[fill] (1.7,4.2) circle (2pt);
		\draw[fill] (2.2,3.8) circle (2pt);
		\draw[fill] (2.7,4.6) circle (2pt);
		\draw[fill] (3.6,3.4) circle (2pt);
		\draw[fill] (4.0,3) circle (2pt);
		\draw[thick, gray] (0.7, 3.6) rectangle (3.4, 5.2);
		\draw[thick, gray] (3.4,3.6) rectangle (4.2, 2.8);
				
		\draw[thick, gray, fill=white] (1.7, 1.8) rectangle (2.2, 2.3);
		\draw[thick, gray, fill=white] (2.2, 2.3) rectangle (2.7, 2.8);
		
		\draw[fill] (1.95, 2.05) circle (2pt);
		\draw[fill] (2.45, 2.55) circle (2pt);

	\end{tikzpicture}
	\quad
	\begin{tikzpicture}[scale=0.9]
		
		\fill[thick, gray, pattern=north west lines, pattern color=lightgray] (0.7, 0.8) rectangle (3.4, 3.6);
		
		\draw[thick, gray, dashed] (1.2, 5) -- (1.2, 0.8);
		\draw[thick, gray, dashed] (1.7, 4.2) -- (1.7, 0.8);
		\draw[thick, gray, dashed] (2.2, 3.8) -- (2.2, 0.8);
		\draw[thick, gray, dashed] (2.7, 4.6) -- (2.7, 0.8);
		
		\draw[fill] (1.2,5) circle (2pt);
		\draw[fill] (1.7,4.2) circle (2pt);
		\draw[fill] (2.2,3.8) circle (2pt);
		\draw[fill] (2.7,4.6) circle (2pt);
		\draw[fill] (3.6,3.4) circle (2pt);
		\draw[fill] (4.0,3) circle (2pt);
		\draw[thick, gray] (0.7, 3.6) rectangle (3.4, 5.2);
		\draw[thick, gray] (3.4,3.6) rectangle (4.2, 2.8);
				
		\draw[thick, gray, fill=white] (2.2, 2.3) rectangle (2.7, 2.8);
		
		\draw[fill] (2.45, 2.55) circle (2pt);

	\end{tikzpicture}
	\caption{}
	\label{figure:Av-2413-3412-skew-new-layer}
\end{center}
\end{figure}

In the leftmost diagram, the initial decreasing run of the new permutation consists entirely of the initial decreasing run of the leftmost component added in the new layer. In the second diagram from the left, the initial decreasing run of the new permutation consists of the initial decreasing run of the new leftmost component together with the leftmost entry of $\pi$. In the rightmost diagram, the initial decreasing run of the new permutation consists of the initial decreasing run of the new leftmost component together with all three entries composing the initial decreasing run of $\pi$. 

We can represent the operation of adding a layer to this particular example by the transformation
\[
	t^3 \longmapsto \left(\frac{z}{1-z}\right)(f(z,t)-1)\left(f(z,1)^3 + tf(z,1)^2 + t^2f(z,1) + t^3f(z,1)\right).
\]
The four summands in the rightmost term represent, when multiplied by $f(z,t)-1$, the four cases of Figure~\ref{figure:Av-2413-3412-skew-new-layer}, from left to right. In each term, $f(z,t)-1$ counts the leftmost nonempty component, the power of $t$ tracks the number of entries from the initial decreasing run of $\pi$ that are part of the initial decreasing run of the new permutation, and the power of $f(z,1)$ tracks the number of components other than the leftmost nonempty component. We have substituted $t=1$ in this last term because the entries in the initial decreasing run of these components are not part of the initial decreasing run of the original permutation. This can be written as a linear operator on monomials:
\begin{align*}
	\Omega : t^k \longmapsto \, & \left(\frac{z}{1-z}\right)(f(z,t)-1)\left(f(z,1)^{k} + tf(z,1)^{k-1} + \cdots + t^{k}\right)\\
	= \, &  \left(\frac{z}{1-z}\right)(f(z,t)-1)\left(\frac{f(z,1)^{k+1}-t^{k+1}}{f(z,1)-t}\right).
\end{align*}

This extends to an operator on power series defined by
\[
	\Omega\left[f^{\not\ominus}(z,t)\right] = \left(\frac{z}{1-z}\right)\left(f(z,t)-1\right)\left(\frac{f(z,1)f^{\not\ominus}(z,f(z,1)) - tf^{\not\ominus}(z,t)}{f(z,1)-t}\right).
\]
As previously noted, every skew indecomposable permutation in $\CC$ that has at least one entry in the region $R_3$ (in the sense of Figure~\ref{figure:Av-2413-3412-skew-1}) can be formed by adding a layer to some smaller skew indecomposable permutation in $\CC$. It follows that $g_{(ii)}(z,t) = \Omega\left[f^{\not\ominus}\right]$, and combined with the equality $f^{\not\ominus}(z,t) = g_{(i)}(z,t) + g_{(ii)}(z,t)$ we have proved the functional equation
\begin{equation}
	\label{equation:fominus-func-1}
	f^{\not\ominus}(z,t) = zf(z,1)(f(z,t)-1) + \left(\frac{z}{1-z}\right)\left(f(z,t)-1\right)\left(\frac{f(z,1)f^{\not\ominus}(z,f(z,1)) - tf^{\not\ominus}(z,t)}{f(z,1)-t}\right).
\end{equation}

We lastly note that Equation~(\ref{equation:f-func-1-no-t}), which relates $f(z)$ and $f^{\not\ominus}(z)$, is easily adapted to track the number of entries in the initial decreasing run. If $\pi \in \CC$ is strictly decreasing, then all entries are part of the initial decreasing run. Otherwise, the initial decreasing run of $\pi$ is composed of the initial decreasing run of the non-trivial skew-indecomposable component as well as all entries of the preceding decreasing permutation. Hence,
\begin{equation}
	\label{equation:f-func-2-with-t}
	f(z,t) = \frac{1}{1-zt} + \frac{f^{\not\ominus}(z,t)}{(1-zt)(1-z)}.
\end{equation}

\subsection{Deriving the generating function}



We employ the ``guess-and-check method'' to prove that $f^{\ominus}(z,t)$ and $f(z,t)$ are algebraic of degree $3$. This application of this method requires several steps.
\begin{enumerate}
	\item Use Equations (\ref{equation:fominus-func-1}) and (\ref{equation:f-func-2-with-t}) to derive a polynomial equation involving the bivariate series $f^{\ominus}(z,t)$, the univariate series $f(z,1)$ and $f^{\ominus}(z,f(z,1))$, and the indeterminate $z$ and $t$.
	\item Guess, with the help of a \texttt{Maple} program, algebraic generating functions (in the form of minimal polynomials) for $f(z,1)$ and $f^{\ominus}(z,f(z,1))$. 
	\item Substitute both guesses into the equation found in Step (1) to find a conjectured minimal polynomial for $f^{\ominus}(z,t)$, and use conjectured this minimal polynomial (and Equation~(\ref{equation:f-func-2-with-t})) to compute minimal polynomials $f(z,1)$ and $f^{\ominus}(z,f(z,1))$. 
	\item Upon verifying that the computed minimal polynomials for $f(z,1)$ and $f^{\ominus}(z,f(z,1))$ are equal to their conjectured guesses and checking several details, the guesses become rigorously verified.
\end{enumerate}

First we inspect Equations~\eqref{equation:fominus-func-1} and \eqref{equation:f-func-2-with-t}. From Equation~\eqref{equation:fominus-func-1} we see that the constant term of $f^{\not\ominus}(z,t)$ is $0$ and hence the constant term of $f(z,t)$ is $1$. Moreover, the form of the equations allows us to determine the coefficient of $z^k$ of both $f(z,t)$ and $f^{\not\ominus}(z,t)$ (which in each case is a polynomial in $t$) from the lower-order terms. This allows us to conclude that there is a \emph{unique} pair of formal power series in $\C[[z,t]]$ that satisfies the pair of equations.

Moreover, one can use this notion of iteration to determine as many initial terms of $f(z,t)$ and $f^{\not\ominus}(z,t)$ as desired (though this quickly becomes computationally cumbersome). Upon finding the coefficients of $z^i$ for $0 \leq i \leq 40$ of $f(z,t)$ and $f^{\not\ominus}(z,t)$, we make the substitutions $t=1$ to the former to find the expansion of $f(z,1)$ up to order $40$, and then we substitute $t=f(z,1)$ into the latter to find the expansion of $f^{\not\ominus}(z,f(z,1))$ up to order $40$.

There are several tools that can be used to automatically conjecture the minimal polynomial for an algebraic generating function given some known initial terms, including \texttt{Gfun}~\cite{salvy:gfun} in \texttt{Maple}, \texttt{Guess}~\cite{kauers:guess} in \texttt{Mathematica}, and (a different) \texttt{Guess}~\cite{rubey:fricas-guess} in \texttt{Fricas}. We instead use \texttt{GuessFunc}, a work-in-progress \texttt{Maple} tool written by the second author. It conjectures nearly instantly minimal polynomials for the two univariate series:
\begin{align}
	z^4\overline{f(z,1)}^3 + (5z^3-11z^2)\overline{f(z,1)}^2 + (3z^2+10z-1)\overline{f(z,1)} - 9z + 1& = 0 \label{equation:min-poly-fz1}\\
	(9z^5-z^4)\overline{f^{\not\ominus}(z,f(z,1))}^3 + (-3z^6+65z^5-168z^4+117z^3-11z^2)\overline{f^{\not\ominus}(z,f(z,1))}^2 \nonumber \\
	+ (-5z^7-27z^6+179z^5-264z^4+116z^3+10z^2-10z+1)\overline{f^{\not\ominus}(z,f(z,1))}\nonumber\\ -z^8-21z^7+66z^6-61z^5+12z^4+6z^3-z^2 & = 0 \label{equation:min-poly-fszfz1}
\end{align}
%

The use of the overline in the equations above denotes that a value is still conjectured. We note now that both of the polynomials have a single solution that is a genuine formal power series in $z$. This can be verified, for example, by examining the structure of the three branches at the origin or by using the trick in Example 6.1.10(d) of Stanley's \emph{Enumerative Combinatorics}~\cite{stanley:ec2}. The conjectured power series solutions $\overline{f(z,1)}$ and $\overline{f^{\not\ominus}(z,f(z,1))}$ imply, upon substituting into Equations~\eqref{equation:fominus-func-1} and \eqref{equation:f-func-2-with-t}, conjectured values of $\overline{f(z,t)}$ and $\overline{f^{\not\ominus}(z,t)}$ that are again formal power series. Because the explicit expressions for $\overline{f(z,1})$ and $\overline{f^{\not\ominus}(z,f(z,1))}$ are unwieldy, we shall use resultants to determine minimal polynomials for $\overline{f(z,t)}$ and $\overline{f^{\not\ominus}(z,t)}$. 

Solving Equation~\eqref{equation:fominus-func-1} for $f^{\not\ominus}(z,t)$ and substituting the result into Equation~\eqref{equation:f-func-2-with-t} gives
\begin{equation*}
	f(z,t) = \frac{1}{1-zt} + \frac{zf(z,1)(f(z,t)-1)((1-z)f(z,1) + f^{\not\ominus}(z,f(z,1)) - t(1-z))}{(ztf(z,t) - (1-z)f(z,1)+t)(1-zt)(-1-z)}.
\end{equation*}
After collecting all terms on one side and clearing the denominator, we find a polynomial
\[
	P_1(y_0, y_1, y_2, x_0, x_1)
\]
such that
\begin{equation*}
	P_1(f(z,t), f(z,1), f^{\not\ominus}(z,f(z,1)), z, t) = 0.
\end{equation*}
As we are seeking to find a minimal polynomial for the $\overline{f(z,t)}$ implied by the conjectures for $\overline{f(z,1)}$ and $\overline{f^{\not\ominus}(z,f(z,1))}$, we also have $P_1(\overline{f(z,t)}, \overline{f(z,1)}, \overline{f^{\not\ominus}(z,f(z,1))}, z,t) = 0$. 

Similarly, the minimal polynomials in Equations~\eqref{equation:min-poly-fz1} and \eqref{equation:min-poly-fszfz1} constitute polynomials $P_2(y_1, x_0)$ and $P_3(y_2, x_0)$, respectively, such that
\[
	P_2(\overline{f(z,1)},z) = 0 \qquad\qquad \text{and} \qquad\qquad P_3(\overline{f^{\not\ominus}(z,f(z,1))}, z) = 0.
\]
By performing elimination via resultants, we are able to find a polynomial $P_4(y_0, x_0, x_1)$ such that \linebreak $P_4(\overline{f(z,t)},z,t)=0$ as follows. First, set $R_1$ to be the resultant of $P_1$ and $P_2$ with respect to $y_1$. We find that $R_1$ factors as
\[
	R_1(y_0, y_2, x_0, x_1) = x_0^2S_1(y_0, y_2, x_0, x_1),
\]
where $S_1$ is irreducible. We thus know that $S_1(\overline{f(z,t)}, \overline{f^{\not\ominus}(z,f(z,1))}, z, t) = 0$. Now set $R_2$ to be the resultant of $S_1$ and $P_3$ with respect to $y_2$. Then, $R_2$ factors as
\[
	R_2(y_0, x_0, x_1) = x_0^{14}(x_0-1)^9(9x_0-1)S_2(y_0, x_0, x_1)S_3(y_0, x_0, x_1),
\]	
where $S_2$ and $S_3$ are irreducible. By computing the first few terms in the series expansions of solutions of $S_2 = 0$ and $S_3 = 0$, we find that the power series expansion of $\overline{f(z,t)}$ is a solution to the equation
\[
	S_2(\overline{f(z,t)},z,t) = 0.
\]
The polynomial $S_2$ has degree $6$ in $y_0$. We can now use $S_2$ to find the minimal polynomial for $\overline{f^{\not\ominus}(z,t)}$. Equation~\eqref{equation:f-func-2-with-t} implies that
\[
	\overline{f^{\not\ominus}(z,t)} - (1-z)(1-zt)\overline{f(z,t)} - z + 1 = 0.
\]
Setting $P_5(y_0, y_3, x_0, x_1) = y_3 - (1-x_0)(1-x_0x_1)y_0 - x_0+1$, the resultant of $S_2$ and $P_5$ with respect to $y_0$ produces a polynomial $S_3(y_3, x_0, x_1)$ such that
\[
	S_3(\overline{f^{\not\ominus}(z,t)}, z, t) = 0.
\]

Let us now summarize the situation. We have in hand minimal polynomials for two bivariate power series $\overline{f(z,t)}$ and $\overline{f^{\not\ominus}(z,t)}$ whose values we conjecture to be equal to $f(z,t)$ and $f^{\not\ominus}(z,t)$. Earlier, we pointed out that there is a unique pair of formal power series solutions to the pair of Equations~\eqref{equation:fominus-func-1} and \eqref{equation:f-func-2-with-t}. So, if we can check that the pair of power series $\overline{f(z,t)}$ and $\overline{f^{\not\ominus}(z,t)}$ is a solution, then they must be the unique solution, and our analysis is complete. In fact, this check is easy to accomplish. Substituting $t=1$ into the minimal polynomial for $\overline{f(z,t)}$ gives
\[
	S_2\left(\left[\overline{f(z,t)}\right]_{t=1}, z, 1\right) = (z-1)^3\left(\left[\overline{f(z,t)}\right]_{t=1}\right)^3P_2\left(\left[\overline{f(z,t)}\right]_{t=1}, z\right) = 0,
\]
proving that $\left[\overline{f(z,t)}\right]_{t=1} = \overline{f(z,1)}$. A similar check (which requires another resultant computation in order to substitute $t=\overline{f(z,1)}$) proves that $\left[\overline{f^{\not\ominus}(z,t)}\right]_{t=f(z,1)} = \overline{f^{\not\ominus}(z,f(z,1))}$. These two verifications complete the ``guess-and-check'' method because they prove that $\overline{f(z,t)}$ are $\overline{f^{\not\ominus}(z,t)}$ are the unique pair of power series solutions to the given functional equations.

Thus we can finally conclude that the generating function for the class $\Av(2413, 3412)$ is $F(z) = f(z,1)$, an algebraic function of degree three whose minimal polynomial is
\[
	z^4F(z)^3 + (5z^3-11z^2)F(z)^2 + (3z^2+10z-1)F(z) - 9z + 1.
\]
We note only in passing that analytic methods (e.g., the ACA Algorithm of Flajolet and Sedgewick~\cite[Section VII.7]{flajolet:ac}) can be used to show that the exponential growth rate of the coefficients of $F(z)$ is exactly $32/5$. A full asymptotic expansion to any desired precision can be similarly calculated.

\section{$\Av(1432,2143)$}
\label{section:Av(1432,2143)}

The left-to-right minima of a permutation are those entries $\pi(\ell)$ such that there does not exist $i < \ell$ with $\pi(i) < \pi(\ell)$. In the plot of a permutation, these are the points for which there is no other point that lies both below and to their left. We call the entries that lie above (and including) the $k$th left-to-right minimum and lie properly below the $(k-1)$th left-to-right minimum the \emph{$k$th slice} of the permutation. Figure~\ref{figure:Av-1432-2143-slices} shows a permutation in $\Av(1432,2143)$ with four slices.

\begin{figure}
	\begin{center}
	\minipage{0.5\textwidth}
		\begin{center}
		\begin{tikzpicture}[scale=0.33]
			\draw[fill] (1,11) circle (6pt);
			\draw[fill] (2,14) circle (6pt);
			\draw[fill] (3,6) circle (6pt);
			\draw[fill] (4,7) circle (6pt);
			\draw[fill] (5,10) circle (6pt);
			\draw[fill] (6,8) circle (6pt);
			\draw[fill] (7,12) circle (6pt);
			\draw[fill] (8,2) circle (6pt);
			\draw[fill] (9,5) circle (6pt);
			\draw[fill] (10,3) circle (6pt);
			\draw[fill] (11,9) circle (6pt);
			\draw[fill] (12,4) circle (6pt);
			\draw[fill] (13,13) circle (6pt);
			\draw[fill] (14,1) circle (6pt);
			\draw[thick, gray] (0.5,0.5) rectangle (14.5,14.5);
			\draw[gray] (0.5,10.5) -- (14.5,10.5);
			\draw[gray] (0.5,5.5) -- (14.5,5.5);
			\draw[gray] (0.5,1.5) -- (14.5,1.5);
			\node at (11, 11) {\footnotesize \emph{$1$st slice}};
			\node at (11, 6) {\footnotesize \emph{$2$nd slice}};
			\node at (11, 2) {\footnotesize \emph{$3$rd slice}};
			\node at (11, 1) {\footnotesize \emph{$4$th slice}};
		\end{tikzpicture}
		\end{center}
		\caption{}
		\label{figure:Av-1432-2143-slices}
	\endminipage\hfill
	\minipage{0.5\textwidth}
		\begin{center}
		\begin{tikzpicture}[scale=.77]
			\draw[thick, lightgray] (0,2) -- (6,2);
			\draw[thick, lightgray] (0,4) -- (6,4);
			\draw[thick, lightgray] (2,0) -- (2,6);
			\draw[thick, lightgray] (4,0) -- (4,6);
			\fill[pattern=north west lines, pattern color=lightgray] (0,4)--(2,4)--(2,2)--(4,2)--(4,0)--(6,0)--(6,4)--(4,4)--(4,6)--(0,6)--cycle;
			
			\draw[thick, gray] (0,0) rectangle (6,6);
			
			\draw[fill] (2,4) circle (2pt) node[above right] {$\pi(a)$};
			\draw[fill] (4,2) circle (2pt) node[above left] {$\pi(b)$};
			\draw[thick] (0.2, 0.2) -- (1.8, 1.8);
			\draw[thick] (0.2, 2.2) -- (1.8, 3.8);
			\draw[thick] (2.2, 0.2) -- (3.8, 1.8);
			\draw[thick] (4.2, 4.2) -- (5.8, 5.8);
			
			\node[above left] at (2, 0) {$R_1$};
			\node[above left] at (4, 0) {$R_2$};
			
		\end{tikzpicture}
		\end{center}
		\caption{}
		\label{figure:Av-1432-2143-one-slice-1}
	\endminipage
	\end{center}
\end{figure}

In order to derive a functional equation for the generating function of $\Av(1432, 2143)$, we will describe 
\begin{enumerate}[(i)]
	\item the permutations consisting of a single slice, i.e., those that start with their smallest entry,
	\item how a slice can be added to an existing permutation.
\end{enumerate}
Noting that every permutation in a class $\CC$ either starts with its smallest entry or can be constructed by adding a slice to a non-empty permutation in $\CC$, this structural description suffices to produce a functional equation.

\subsection{Single-slice permutations}
\label{subsection:single-slice}

A permutation in $\Av(1432, 2143)$ that starts with its smallest entry has the form $\pi = 1 \oplus \tau$ for $\tau \in \Av(321, 2143)$. With the benefit of foresight, we closely examine three separate cases of permutations that have this form: $\tau$ is increasing, $\tau$ contains a $21$ pattern but not a $2413$ pattern, and $\tau$ contains a $2413$ pattern.

The first case is clear. For the second case, suppose now that $\tau$ contains a $21$ pattern but not a $2413$ pattern. Let $\pi(a)$ be the largest entry that plays the role of the $2$ in a $21$ pattern, and let $\pi(b)$ be the rightmost entry such that $\pi(a)\pi(b)$ forms a $21$ pattern.  Figure~\ref{figure:Av-1432-2143-one-slice-1} shows the plot of such a permutation. Each of the non-shaded regions must hold an increasing permutation; otherwise a $321$ or $2143$ pattern is created.

Moreover, because $\tau$ avoids $2413$ there cannot be a $21$ pattern for which the $2$ is in the region $R_1$ and the $1$ is in the region $R_2$. Figure~\ref{figure:Av-1432-2143-one-slice-2} shows the refined diagram. One can clearly see that any permutation drawn on this diagram will avoid the patterns $321$, $2143$, and $2413$, completing the description of this case.

\begin{figure}
\begin{center}
	\minipage{0.5\textwidth}
		\begin{center}
		\begin{tikzpicture}[scale=0.77]
			
			\draw[thick, lightgray] (0,2) -- (6,2);
			\draw[thick, lightgray] (0,4) -- (6,4);
			\draw[thick, lightgray] (2,0) -- (2,6);
			\draw[thick, lightgray] (4,0) -- (4,6);
			\fill[pattern=north west lines, pattern color=lightgray] (0,4)--(2,4)--(2,2)--(4,2)--(4,1)--(2,1)--(2,2)--(0,2)--(0,1)--(2,1)--(2,0)--(6,0)--(6,4)--(4,4)--(4,6)--(0,6)--cycle;
			
			\draw[thick, lightgray, dashed] (0,1)--(6,1);
			\draw[thick, gray] (0,0) rectangle (6,6);
			
			\draw[fill] (2,4) circle (2pt) node[above right] {$\pi(a)$};
			\draw[fill] (4,2) circle (2pt) node[above left] {$\pi(b)$};
			\draw[thick] (0.2, 0.2) -- (1.8, 0.8);
			\draw[thick] (0.2, 2.2) -- (1.8, 3.8);
			\draw[thick] (2.2, 1.2) -- (3.8, 1.8);
			\draw[thick] (4.2, 4.2) -- (5.8, 5.8);
		\end{tikzpicture}
		\caption{}
		\label{figure:Av-1432-2143-one-slice-2}
		\end{center}
	\endminipage\hfill
	\minipage{0.5\textwidth}
		\begin{center}
		\begin{tikzpicture}[scale=0.77]
			\fill[pattern=north west lines, pattern color=lightgray] (0,1.5)--(1.5,1.5)--(1.5,0)--(7.5,0)--(7.5,6)--(6,6)--(6,7.5)--(0,7.5)--cycle;
			\fill[white] (1.5,3) rectangle (2.25,4.5);
			\fill[white] (2.25,4.5) rectangle (3,6);
			\fill[white] (4.5,1.5) rectangle (5.25, 3);
			\fill[white] (5.25, 3) rectangle (6,4.5);
					
			\foreach \i in {1,2,3,4}{
				\draw[thick, lightgray] (0,{1.5*\i}) -- (7.5,{1.5*\i});
			}
			\foreach \j in {1,2,3,4}{
				\draw[thick, lightgray] ({1.5*\j},0) -- ({1.5*\j},7.5);
			}
			
			\draw[thick, lightgray, dashed] (2.25,0)--(2.25,7.5);
			\draw[thick, lightgray, dashed] (5.25,0)--(5.25,7.5);
			
			\draw[thick, gray] (0,0) rectangle (7.5,7.5);
			\draw[fill] (1.5,3) circle (2pt);
			\draw[fill] (3,6) circle (2pt);
			\draw[fill] (4.5,1.5) circle (2pt);
			\draw[fill] (6,4.5) circle (2pt);

			\draw[thick] (0.2, 0.2) -- (1.3, 1.3);
			\draw[thick] (6.2, 6.2) -- (7.3, 7.3);
			\draw[thick] (6.2, 6.2) -- (7.3, 7.3);
			\draw[thick] (1.7, 3.2) -- (2.05, 4.3);
			\draw[thick] (2.45, 4.7) -- (2.8, 5.8);
			\draw[thick] (4.7, 1.7) -- (5.05, 2.8);
			\draw[thick] (5.45, 3.2) -- (5.8, 4.3);		
		
	\end{tikzpicture}
	\caption{}
	\label{figure:Av-1432-2143-one-slice-3}
		\end{center}
	\endminipage
\end{center}
\end{figure}

The last case to consider is when $\tau$ contains a $2413$ pattern. Figure~\ref{figure:Av-1432-2143-one-slice-3} shows the diagram for such a permutation, where the $4$ is chosen to be as large as possible, the $1$ as small as possible for the given $4$, the $2$ as leftmost as possible for the given $1$ and $4$, and the $3$ as rightmost as possible for the given $1$, $2$, and $4$. Any permutation drawn on this diagram avoids $321$ and $2143$ and so this completes the classification of the permutations in $\Av(1432, 2143)$ that contain only one slice.

\subsection{Multi-slice permutations}

A permutation in $\Av(1432,2143)$ consisting of multiple slices can be built by adding a single slice of the form in the previous subsection to a permutation in $\Av(1432,2143)$ with one fewer slices. The existing entries in the smaller permutation restrict where the entries of the new slice can lie, and different permutations can be formed from the same new slice by interleaving the entries in a different way.

The \emph{trailing increasing sequence} of a permutation is the longest sequence of consecutive increasing entries at the rightmost end of the permutation. The \emph{gaps} of a permutation are the vertical regions that the lay to the right of the rightmost $2$ in any $21$ pattern and in between entries of the trailing increasing sequence. A permutation with a trailing increasing sequence of length $\ell$ has $\ell+1$ gaps. (Note that every permutation has at least two gaps.)

To prevent the creation of a $2143$ pattern, any entries in an added slice must lie within the gaps of the permutation to which they are being added. We consider three different ways to add a slice to a permutation in $\Av(1432, 2143)$.
\begin{itemize}
	\item \emph{Case 1:} All added entries are contained in a single gap.
	\item \emph{Case 2:} At least two gaps are used, and there is a $1324$ pattern in the added slice such that the $1$ and the $3$ are contained in one gap and the $2$ and the $4$ are contained in another gap.
	\item \emph{Case 3:} At least two gaps are used, and there is no such $1324$ pattern. 
\end{itemize}

Case 1 is the easiest. To add such a slice to a permutation $\pi$, one just needs to choose which of the gaps to insert it into. 

To insert a new slice as in Case 2 we must first choose two gaps and then insert four entries that form a $1324$ such that the $1$ and $3$ are in one gap while the $2$ and $4$ are in another. This leads to a diagram like that in Figure~\ref{figure:Av-1432-2143-multi-slice-1}, in which the first slice is the permutation $14235$ with trailing increasing sequence of length $3$ and hence $4$ gaps. The $1$ in the $1324$ is the first entry of the new slice (the newly added left-to-right minimum), the $3$ is chosen to be as high as possible for the given $1$, then $2$ as low as possible, and the $4$ as high as possible. Note that we have chosen, in this example, two non-consecutive gaps. Once such a $1324$ pattern is added, no entries may be inserted into the other gaps without creating a $1432$ or $2143$ pattern. 

\begin{figure}
\begin{center}
	\begin{tikzpicture}[scale=1]

		\fill[pattern=north west lines, pattern color=lightgray] (0,1)--(7,1)--(7,6)--(0,6)--cycle;
		\fill[white] ({1+2/3},2) rectangle ({1+4/3},4);
		\fill[white] ({1+4/3},5) rectangle (3,6);
		\fill[white] ({4+2/3},3) rectangle ({4+4/3},5);
		
				
		

		\draw[thick, lightgray] (0,6)--(7,6);
		\draw[thick, lightgray] (1,6.75)--(1,1);
		\draw[thick, lightgray] (3,6.25)--(3,1);
		\draw[thick, lightgray] (4,6.5)--(4,1);
		\draw[thick, lightgray] (6,7)--(6,1);
		
		\draw[thick, lightgray] (0,2)--(7,2);
		\draw[thick, lightgray] (0,3)--(7,3);
		\draw[thick, lightgray] (0,4)--(7,4);
		\draw[thick, lightgray] (0,5)--(7,5);
		
		\draw[thick, lightgray] ({1+2/3},1)--({1+2/3},6);
		\draw[thick, lightgray] ({1+4/3},1)--({1+4/3},6);
		\draw[thick, lightgray] ({4+2/3},1)--({4+2/3},6);
		\draw[thick, lightgray] ({4+4/3},1)--({4+4/3},6);

		\draw[thick, gray] (0,1) rectangle (7,7);
		\draw[fill] (0,6) circle (2pt);
		\draw[fill] (1,6.75) circle (2pt);
		\draw[fill] (3,6.25) circle (2pt);
		\draw[fill] (4,6.5) circle (2pt);
		\draw[fill] (6,7) circle (2pt);
		
		\draw[fill] ({1+2/3},2) circle (2pt);
		\draw[fill] ({1+4/3},4) circle (2pt);
		\draw[fill] ({4+2/3},3) circle (2pt);
		\draw[fill] ({4+4/3},5) circle (2pt);

		\node[below=2pt, above] at (2,6) {\footnotesize gap $1$};
		\node[below=2pt, above] at (3.5,6) {\footnotesize gap $2$};
		\node[below=2pt, above] at (5,6) {\footnotesize gap $3$};
		\node[below=2pt, above] at (6.5,6) {\footnotesize gap $4$};

		\draw[thick] ({1+2/3+0.2},2.2) -- ({1+4/3-0.2},3.8);
		\draw[thick] ({4+2/3+0.2},3.2) -- ({4+4/3-0.2},4.8);
		\draw[thick] ({1+4/3+0.2},5.2) -- (2.8,5.8);
		
	\end{tikzpicture}
	\caption{}
	\label{figure:Av-1432-2143-multi-slice-1}
\end{center}
\end{figure}

All permutations that can be drawn on Figure~\ref{figure:Av-1432-2143-multi-slice-1} avoid the patterns $1432$ and $2143$ and therefore the classification of Case 2 is complete.

In Case 3, the newly added slice must use at least two gaps and must not contain the kind of $1324$ pattern described in Case 2. This forces the added slice to contain only two increasing sequences stacked on top of each other---the lower one can span any number of gaps, while the upper one can span only the leftmost used gap. Figure~\ref{figure:Av-1432-2143-multi-slice-2} demonstrates the diagram of such a permutation. It assumes that the leftmost used gap is the second gap, though this need not be the case. Once again, all permutations that can be drawn on this diagram avoid the patterns $1432$ and $2143$, and hence the analysis of the three cases is complete.

\begin{figure}
\begin{center}
	\begin{tikzpicture}[scale=1]

		\fill[pattern=north west lines, pattern color=lightgray] (0,1)--(7,1)--(7,6)--(0,6)--cycle;
		\fill[white] ({2+2/3},2) rectangle (4,3);
		\fill[white] (4,3) rectangle (5.5,4);
		\fill[white] (5.5,4) rectangle (7,5);
		\fill[white] ({2+2/3},5) rectangle (4,6);
		
		\draw[thick, lightgray] (0,6)--(7,6);
		\draw[thick, lightgray] (1,6.75)--(1,1);
		\draw[thick, lightgray] (2,6.25)--(2,1);
		\draw[thick, lightgray] (4,6.5)--(4,1);
		\draw[thick, lightgray] (5.5,7)--(5.5,1);
		
		\draw[thick, lightgray] (0,2)--(7,2);
		\draw[thick, lightgray] (0,3)--(7,3);
		\draw[thick, lightgray] (0,4)--(7,4);
		\draw[thick, lightgray] (0,5)--(7,5);
		
		\draw[thick, lightgray] ({2+2/3},1)--({2+2/3},6);
		
		\draw[thick, gray] (0,1) rectangle (7,7);
		\draw[fill] (0,6) circle (2pt);
		\draw[fill] (1,6.75) circle (2pt);
		\draw[fill] (2,6.25) circle (2pt);
		\draw[fill] (4,6.5) circle (2pt);
		\draw[fill] (5.5,7) circle (2pt);
		
		\draw[fill] ({2+2/3},2) circle (2pt);
		\draw[fill] (4.2,3.2) circle (2pt);
		
		\node[below=2pt, above] at (1.5,6) {\footnotesize gap $1$};
		\node[below=2pt, above] at (3,6) {\footnotesize gap $2$};
		\node[below=2pt, above] at (4.75,6) {\footnotesize gap $3$};
		\node[below=2pt, above] at (6.25,6) {\footnotesize gap $4$};
		
		\draw[thick] ({2+2/3+0.2},2.2)--(3.8, 2.8);
		\draw[thick] (4.2, 3.2) -- (5.3, 3.8);
		\draw[thick] (5.7, 4.2) -- (6.8, 4.8);
		\draw[thick] ({2+2/3+0.2},5.2)--(3.8, 5.8);
		
	\end{tikzpicture}
	\caption{}
	\label{figure:Av-1432-2143-multi-slice-2}
\end{center}
\end{figure}

\subsection{Constructing a functional equation}

As observed above, in order to add a new slice to a permutation $\pi$ one must know the number of gaps that $\pi$ has. To that end, we define $f(z,t)$ to be the bivariate generating function for permutations in $\Av(1432, 2143)$ where $z$ tracks the length of $\pi$ and $t$ tracks one fewer than the number of gaps in $\pi$---equivalently, $t$ tracks the length of the trailing increasing sequence. 

In this section we will derive a functional equation for $f(z,t)$. We noted earlier that every permutation in a class is either a single-slice permutation, or can be built by adding a slice to a permutation in the same class with one fewer slices. This can be simplified if one considers the empty permutation of length zero to be a single-slice permutation with a single gap. Now, every nonempty permutation is the result of adding a slice to either the empty permutation or a permutation with one fewer slice. The cases considered in the previous section lead to the natural split
\[
	f(z,t) = 1 + g_{(a)}(z,t) + g_{(b)}(z,t) + g_{(c)}(z,t),
\]
where $g_{(a)}$, $g_{(b)}$, and $g_{(c)}$ are respectively the generating function for permutations built as in Cases 1, 2, and 3 above.

Before finding formulas for these three functions, we need to find the generating function for the single-slice permutations. It is important that each entry in the trailing increasing sequence is marked by a power of $t$. The simplest single-slice permutation is a strictly increasing permutation. These have generating function
\[
	\frac{z}{1-tz}.
\]
In this one special case, we do not count the leftmost entry as part of the trailing increasing sequence because there cannot be a gap to the left of the first entry as this would invalidate the assumption that the leftmost entry is a left-to-right minimum.

The single-slice permutation $1 \oplus \tau$ where $\tau$ has the form shown in Figure~\ref{figure:Av-1432-2143-one-slice-2} have generating function
\[
	tz^3 \left(\frac{1}{1-2z}\right)\left(\frac{1}{1-tz}\right)^2.
\]
The $tz^3$ term counts the three placed points (the leftmost point of $\pi$, and both points that form the $21$ pattern), only one of which is part of the trailing increasing sequence. The $1/(1-2z)$ term counts the number of ways to draw a permutation in the two cells in the leftmost column of Figure~\ref{figure:Av-1432-2143-one-slice-2}. Lastly, the $1/(1-tz)^2$ term counts the number of ways to draw entries on the increasing line segments in columns 2 and 3, all of whose entries will be part of the trailing increasing sequence.

Lastly, the generating function for single-slice permutations of the form shown in Figure~\ref{figure:Av-1432-2143-one-slice-3} is
\[
	\frac{t^2z^5}{(1-z)^2(1-tz)^2(1-(1+t)z)}.
\]
Summarizing, the generating function for single-slice permutations, not including the empty permutation, is 
\[
	s(z,t) = \frac{z}{1-tz} + \frac{tz^3}{(1-2z)(1-tz)^2} + \frac{t^2z^5}{(1-z)^2(1-tz)^2(1-(1+t)z)}.
\]	

With this in hand, we derive $g_{(a)}$, $g_{(b)}$, and $g_{(c)}$. The permutations counted by $g_{(a)}$ are those formed by starting with a permutation and adding a slice with entries in only one gap. When starting with a permutation with a trailing increasing sequence of length $k$, there are $k+1$ gap options. Moreover, any entries of the trailing increasing sequence to the right of the chosen gap remain in the trailing increasing sequence after insertion of the new slice. This action can be viewed as a linear operator acting on monomials
\[
	\Phi : t^k \mapsto s(z,t)\left(1 + t + t^2 + \cdots + t^k\right).
\]
On the level of generating functions, this is accomplished by defining
\[
	g_{(a)}(z,t) = \Phi[f(z,t)] = s(z,t)\left(\frac{f(z,1)-tf(z,t)}{1-t}\right).
\]

The derivation of $g_{(b)}$ requires more care. We break down the process of adding a slice containing a $1324$ (split in the middle by an entry from the original permutation) into two parts: first any gap other than the rightmost is chosen for the leftmost entry of the new slice, then a second gap to its right is chosen and all entries are placed. The first step is modeled simply by the linear operator
\[
	\Theta: t^k \mapsto t + t^2 + \cdots + t^k
\]
(with $t^0 \mapsto 0$) whose action on functions is
\[
	\Theta[A(z,t)] = \frac{t(A(z,1)-A(z,t))}{1-t}.
\]
This operator does not actually place the leftmost entry of the new slice. It only picks the gap it will eventually go in, and ``unmarks'' all entries of the trailing increasing sequence of the original permutation to its left.

In the second step, another gap is chosen and a permutation is added into the two gaps as in Figure~\ref{figure:Av-1432-2143-multi-slice-1}. This corresponds to the operator
\[
	\Psi : t^k \mapsto \frac{t^2z^4}{(1-z)^2(1-tz)(1-(1+t)z)}(1 + t + t^2 + \cdots + t^{k-1}).	
\]
The $1 + t + t^2 + \cdots + t^{k-1}$ represents the fact that, based on the choice of second gap, between $0$ and $k-1$ entries from the trailing increasing sequence of the original will remain in the trailing increasing sequence of the new permutation. The action of $\Psi$ on functions is
\[
	\Psi[B(z,t)] =  \frac{t^2z^4}{(1-z)^2(1-tz)(1-(1+t)z)}\left(\frac{B(z,1)-B(z,t)}{1-t}\right).
\]
It follows from our setup that
\[
	g_{(b)}(z,t) = \Psi[\Theta[f(z,t)]].
\]
There is a slight wrinkle in this composition: $\Theta[f(z,t)]$ doesn't appear at first to be defined at $t=1$. In fact, this is the \emph{discrete derivative}, and has formal power series:
\[
	\left[\frac{B(z,1)-B(z,t)}{1-t}\right]_{t=1} = B_t(z,1),
\]
where $B_t$ denotes the partial derivative of $B$ with respect to $t$. This allows us to compute the composition fully:
\[
	g_{(b)}(z,t) = \frac{t^2z^4}{(1-z)^2(1-tz)(1-(1+t)z)}\left(\frac{f_t(z,1)}{1-t} - \frac{t(f(z,1)-f(z,t))}{(1-t)^2}\right).
\]

The computation of $g_{(c)}$ follows along the same lines. First, the gap that will contain the leftmost entry of the new slice is chosen, and all entries are added to it. This is modeled by the linear operator
\[
	\Lambda : t^k \mapsto \frac{z}{1-2z}(t + t^2 + \cdots + t^k),
\]
yielding
\[
	\Lambda[C(z,t)] = \frac{z}{1-2z}\left(\frac{t(C(z,1)-C(z,t))}{1-t}\right).
\]
The second step involves picking a rightmost gap that will contain entries, placing a nonempty increasing sequence of entries (all of which will be part of the new trailing increasing sequence), and adding arbitrary length increasing sequences in all gaps in between (see Figure~\ref{figure:Av-1432-2143-multi-slice-2}). This step produces the linear operator
\begin{align*}
	\Xi : t^k \mapsto &\frac{tz}{1-tz}\left(t^{k-1} + \frac{t^{k-2}}{1-z} + \frac{t^{k-3}}{(1-z)^2} + \cdots + \frac{t}{(1-z)^{k-2}} + \frac{1}{(1-z)^{k-1}}\right)\\
	= & \left(\frac{tz}{1-tz}\right)\left(\frac{t^k - \left(\frac{1}{1-z}\right)^k}{t-\frac{1}{1-z}}\right).
\end{align*}
Allow us to explain this operator a bit more. In the parenthesized part of the first line, the $t^{k-1}$ term represents the situation when the second chosen gap is immediately adjacent to the first gap, leaving no middle gaps for increasing sequences. The $t^{k-2}/(1-z)$ term represents the situation when there is one intermediate gap. On the other extreme, the $1/(1-z)^{k-1}$ term represents the situation when the second gap is chosen as far to the right as possible---all original entries in the trailing increasing sequence are no longer part of the new one, and there are $k-1$ intermediate gaps, all of which hold increasing permutations.

The action of $\Xi$ on power series is
\[
	\Xi[D(z,t)] = \frac{tz}{1-tz}\left(\frac{D(z,t) - D\left(z,\frac{1}{1-z}\right)}{t - \frac{1}{1-z}}\right).
\]
The composition can be computed without the problems encountered previously:
\begin{align*}
	g_{(c)}(z,t) &= \Xi[\Lambda[f(z,t)]]\\
	&= \frac{tz^2}{(1-2z)(1-tz)}\left(\frac{t(f(z,1)-f(z,t))}{(1-t)\left(t-\frac{1}{1-z}\right)} - \frac{\frac{1}{1-z}\left(f(z,1) - f\left(z,\frac{1}{1-z}\right)\right)}{\left(1 - \frac{1}{1-z}\right)\left(t-\frac{1}{1-z}\right)}\right)\\
	&= \frac{tz^2}{(1-2z)(1-tz)}\left(\frac{t(z-1)(f(z,1)-f(z,t))}{(1-t)(1-t+tz)} + \frac{(z-1)\left(f(z,1) - f\left(z,\textstyle\frac{1}{1-z}\right)\right)}{z(1-t+tz)} \right).
\end{align*}

Combining the results from this section, we find the functional equation

\begin{align}
\begin{split}
	\label{equation:av-1432-2143-functional}
	f(z,t) &= 1 + \left(\frac{z}{1-tz} + \frac{tz^3}{(1-2z)(1-tz)^2} + \frac{t^2z^5}{(1-z)^2(1-tz)^2(1-(1+t)z)}\right)\left(\frac{f(z,1)-tf(z,t)}{1-t}\right)\\[3pt]
		&\quad + \frac{t^2z^4}{(1-z)^2(1-tz)(1-(1+t)z)}\left(\frac{f_t(z,1)}{1-t} - \frac{t(f(z,1)-f(z,t))}{(1-t)^2}\right)\\
		&\quad + \frac{tz^2}{(1-2z)(1-tz)}\left(\frac{t(z-1)(f(z,1)-f(z,t))}{(1-t)(1-t+tz)} + \frac{(z-1)\left(f(z,1) - f\left(z,\textstyle\frac{1}{1-z}\right)\right)}{z(1-t+tz)} \right)
\end{split}
\end{align}

\subsection{Deriving the generating function}

A minimal polynomial for $f(z,1)$ can be derived from Equation (\ref{equation:av-1432-2143-functional}) with the techniques of Bousquet-M\'elou and Jehanne~\cite{bousquet-melou:poly-eqs}. By rearranging Equation  (\ref{equation:av-1432-2143-functional}) and taking only the numerator, we find a polynomial $P(y_0, y_1, y_2, y_3, x_0, x_1)$ such that
\begin{equation}
	\label{equation:poly-eq}
	P\left(f(z,t), f(z,1), f_t(z,1), f\left(z, \frac{1}{1-z}\right), z, t\right) = 0.
\end{equation}
To our advantage, $P$ is linear in each of $y_0$, $y_1$, $y_2$, and $y_3$, and we can write
\[
	P(y_0, y_1, y_2, y_3, x_0, x_1) = K(x_0,x_1)y_0 + R(y_1, y_2, y_3, x_0, x_1)
\]
where
\begin{align*}
	K(z,t) &= (1-2z)(1-z)(1-z-t+t^2z-t^2z^2)\\
		&\qquad\qquad (1-z-2t+2tz+t^2+t^2z-3t^2z^2-2t^3z+2t^2z^3+2t^3z^2+t^4z^2-t^4z^3).
\end{align*}
Theorem 2 of~\cite{bousquet-melou:poly-eqs} guarantees that there are exactly $3$ fractional power series $t(z)$ such that $K(z,t(z)) = 0$, and by checking initial terms we find that they are distinct. We call these $t_1(z)$, $t_2(z)$, and $t_3(z)$. More precisely, $t_1(z)$ has minimal polynomial 
\[
	m_1(z,t) = 1-z-t+t^2z-t^2z^2
\]
and $t_2(z)$ and $t_3(z)$ have minimal polynomial
\[
	m_2(z,t) = 1-z-2t+2tz+t^2+t^2z-3t^2z^2-2t^3z+2t^2z^3+2t^3z^2+t^4z^2-t^4z^3.
\]

Substituting $t = t_i(z)$ for $i=1,2,3$ into Equation~\eqref{equation:poly-eq} yields three equations
\begin{align*}
	R\left(f(z,1),f_t(z,1),f\left(z,\frac{1}{1-z}\right), z, t_1(z)\right) &= 0,\\
	R\left(f(z,1),f_t(z,1),f\left(z,\frac{1}{1-z}\right), z, t_2(z)\right) &= 0,\\
	R\left(f(z,1),f_t(z,1),f\left(z,\frac{1}{1-z}\right), z, t_3(z)\right) &= 0.
\end{align*}
These three equations, together with $m_1(z,t_1) = 0$, $m_2(z,t_2) = 0$, and $m_2(z,t_3) = 0$, constitute a system of six equations with six unknowns. At this point one would typically use resultants to eliminate the variables $t_1, t_2, t_3, f_t(z,1), f(z, 1/(1-z))$ from the system, as elimination via Gr\"obner bases tends to take too long. To our surprise, \texttt{Maple 2017} is able to perform the necessary computation in under an hour, revealing that $f(z,1)$ has a minimal polynomial of degree $8$:
\begin{footnotesize}
\begin{align*}
	& \left( 64{z}^{12}-496{z}^{11}+1728{z}^{10}-3672{z}^{9}+5356{z}^{8}-5591{z}^{7}+4244{z}^{6}-2363{z}^{5}+950{z}^{4}-259{z}^{3}+42{z}^{2}-3z \right)  \left( f \left( z,1 \right)  \right) ^{8}\\
	& + \left( -128{z}^{14}+1408{z}^{13}-7040{z}^{12}+21456{z}^{11}-45000{z}^{10}+69088{z}^{9}-79916{z}^{8}+70289{z}^{7}\right.\\
	&\qquad\qquad \left.-46994{z}^{6}+23573{z}^{5}-8438{z}^{4}+1903{z}^{3}-189{z}^{2}-11z+3 \right)  \left( f \left( z,1 \right)  \right) ^{7}\\
	& + \left( 64{z}^{16}-896{z}^{15}+5968{z}^{14}-25040{z}^{13}+73708{z}^{12}-160832{z}^{11}+268543{z}^{10}-348733{z}^{9}\right.\\
	&\qquad\qquad \left.+353476{z}^{8}-278054{z}^{7}+166834{z}^{6}-72514{z}^{5}+19397{z}^{4}-1123{z}^{3}-1132{z}^{2}+356z-34 \right)  \left( f \left( z,1 \right)  \right) ^{6}\\
	& + \left( -128{z}^{16}+1824{z}^{15}-12048{z}^{14}+49120{z}^{13}-138784{z}^{12}+288562{z}^{11}-455345{z}^{10}+550533{z}^{9}\right.\\
	&\qquad\qquad \left.-505169{z}^{8}+338901{z}^{7}-145627{z}^{6}+12147{z}^{5}+35775{z}^{4}-28609{z}^{3}+10771{z}^{2}-2072z+162 \right)  \left( f \left( z,1 \right)  \right)^{5}\\
	& + \left( 16{z}^{17}-176{z}^{16}+728{z}^{15}-816{z}^{14}-4623{z}^{13}+24411{z}^{12}-57292{z}^{11}+71304{z}^{10}-12983{z}^{9}\right.\\
	&\qquad\qquad \left.-134003{z}^{8}+305614{z}^{7}-403001{z}^{6}+374932{z}^{5}-251513{z}^{4}+116646{z}^{3}-34673{z}^{2}+5843z-420 \right)  \left( f \left( z,1 \right)  \right)^{4}\\
	& + \left( -80{z}^{15}+1128{z}^{14}-7732{z}^{13}+34582{z}^{12}-113186{z}^{11}+284672{z}^{10}-560167{z}^{9}+869521{z}^{8}\right.\\
	&\qquad\qquad \left.-1074044{z}^{7}+1061462{z}^{6}-829369{z}^{5}+492560{z}^{4}-208617{z}^{3}+57889{z}^{2}-9257z+639 \right)  \left( f \left( z,1 \right)  \right) ^{3}\\
	& + \left( 8{z}^{16}-96{z}^{15}+470{z}^{14}-598{z}^{13}-6029{z}^{12}+43640{z}^{11}-157743{z}^{10}+380749{z}^{9}-671697{z}^{8}\right.\\
	&\qquad\qquad \left.+906244{z}^{7}-952092{z}^{6}+770540{z}^{5}-462861{z}^{4}+195020{z}^{3}-53381{z}^{2}+8396z-570 \right)  \left( f \left( z,1 \right)  \right) ^{2}\\
	& + \left( -44{z}^{14}+428{z}^{13}-1385{z}^{12}-485{z}^{11}+18705{z}^{10}-74725{z}^{9}+177105{z}^{8}-296154{z}^{7}+365009{z}^{6}\right.\\
	&\qquad\qquad \left.-329611{z}^{5}+211661{z}^{4}-92435{z}^{3}+25723{z}^{2}-4068z+276 \right) f \left(z,1 \right)\\¡
	& +{z}^{15}-5{z}^{14}-10{z}^{13}+130{z}^{12}-378{z}^{11}+230{z}^{10}+2328{z}^{9}-11440{z}^{8}+29440{z}^{7}-48242{z}^{6}+52233{z}^{5}\\
	&\qquad\qquad -37467{z}^{4}+17470{z}^{3}-5050{z}^{2}+816z-56	.
\end{align*}
\end{footnotesize}

Analytic methods can be used to show that the exponential growth rate of $\Av(1432,2143)$ is a root of the polynomial $z^4 - 7z^3 + 9z^2 - 8z + 4$, and is approximately $5.63$.


\section{Concluding Remarks}
\label{section:final-observations}

	For both \twobyfour{} classes that we've enumerated here, the first step was to find a structural description that shed light on how larger permutations in the class may be built from smaller ones. The structural descriptions were then translated to quite complicated functional equations. At this point the methodologies diverge. In Section~\ref{section:Av(2413,3412)}, the functional equation is solved using the guess-and-check paradigm, while in Section~\ref{section:Av(1432,2143)} we use the techniques of~\cite{bousquet-melou:poly-eqs}. This is by necessity. Although the guess-and-check paradigm could theoretically be applied to the functional equation in Section~\ref{section:Av(1432,2143)}, the computations required are too large to be completed in a reasonable length of time.\footnote{To be more precise, we have no problem using the functional equation to produce 250 terms, nor using the 250 terms to automatically conjecture the minimal polynomial for $f(z,1)$. Rather, the ``check'' step is where the computations become unfeasible.} On the other hand, the functional equations in Section~\ref{section:Av(2413,3412)} seem to fail the genericity conditions required to apply the methods from~\cite{bousquet-melou:poly-eqs}. As there are two unknown univariate functions in Equations~\eqref{equation:fominus-func-1} and~\eqref{equation:f-func-2-with-t}, we would hope to have a kernel $K(z,t)$ (akin to the one used in Section~\ref{section:Av(1432,2143)}) for which there are two fractional power series $t_1(z)$ and $t_2(z)$ that cancel the kernel. Instead we find only one. It is possible that there is some manipulation that can be done to these equations to fix this deficiency.

	A listing of all \twobyfour{} classes can be found on \href{https://en.wikipedia.org/wiki/Enumerations_of_specific_permutation_classes}{Wikipedia}~\cite{wiki:enum}. With the two enumerations computed here, the \twobyfour{} classes are almost completely understood. Polynomial-time counting algorithms are now known for all $38$ distinct enumerations. For all except three, the exact generating function is known, and hence complete asymptotic expansions to arbitrary precision can be computed. The remaining three classes---$\Av(4123,4231)$, $\Av(4123,4312)$, and $\Av(4231,4321)$---can be described using a sorting-machine model, as explained by Albert, Homberger, Pantone, Shar, and Vatter~\cite{albert:C-machines}, that yields polynomial-time counting algorithms for all three. Still, their generating functions are not known, and their asymptotic behavior has only been non-rigorously estimated. 
	
	This may be for good reason. Conjecture 5.4 of~\cite{albert:C-machines} posits that the generating function for these three classes may not be D-finite; that is, their generating functions may not be the solution to a linear differential equation with polynomial coefficients. The class of D-finite functions includes as a subset all algebraic functions. If correct, this represents an interesting dichotomy: all \twobyfour{} classes are either algebraic or non-D-finite\footnote{In fact, Conjecture 5.4 goes even further, suggesting the three generating functions don't even belong to the much broader family of D-algebraic functions.}.
	
\textbf{Acknowledgments.} We thank Michael Albert for his open-source software \texttt{PermLab}~\cite{albert:permlab-github}, which was a useful aide in this research.

\bibliographystyle{acm}
\bibliography{/Users/jay/Dropbox/Research/refs/bib.bib}

\end{document}